\documentclass[a4paper,12pt]{article}
\usepackage{amssymb}
\usepackage{verbatim}
\usepackage{amsmath}
\usepackage{amscd}
\usepackage{array}
\usepackage[arrow,matrix]{xy}
\usepackage[bookmarksnumbered,colorlinks,linkcolor=black,citecolor=black,urlcolor=black]{hyperref}

\newtheorem{theo}{Theorem}[section]
\newtheorem{prop}[theo]{Proposition}
\newtheorem{lem}[theo]{Lemma}
\newtheorem{cor}[theo]{Corollary}
\newtheorem{defi}[theo]{Definition}
\newtheorem{rem}[theo]{Remark}
\newtheorem{ex}[theo]{Example}
\newtheorem{con}[theo]{Convention}

\def \Br {{\rm{Br}}}

\def \si {{\sigma}}

\def \Gal {{\rm{Gal}}}

\def \Im {{\rm {Im}}}

\def \Spec {{\rm{Spec}}}

\def \Hom {{\rm {Hom}}}
\def \End {{\rm {End}}}

\def \Aut{{\rm Aut}}

\def\ov{\overline}

\def \Z {{\mathbb Z}}
\def \Q {{\mathbb Q}}
\def \F {{\mathbb F}}

\def \EE {{\rm E}}
\def \UU {{\rm U}}

\def \rk {{\rm{rk}}}

\def\C{{\mathbb C}}

\def\lra{\longrightarrow}

\def\id{{\rm id}}

\def\H{{\rm H}}

\def\Kum{{\rm Kum}}

\def\NS{{\rm NS\,}}

\def\O{{\cal O}}

\def\si{\sigma}

\def\discr{{\rm discr}}

\def\e{\varepsilon}

\newcommand{\bthe}{\begin{theo}}
\newcommand{\ble}{\begin{lem}}
\newcommand{\bpr}{\begin{prop}}
\newcommand{\bco}{\begin{cor}}
\newcommand{\bde}{\begin{defi}}
\newcommand{\ethe}{\end{theo}}
\newcommand{\ele}{\end{lem}}
\newcommand{\epr}{\end{prop}}
\newcommand{\eco}{\end{cor}}
\newcommand{\ede}{\end{defi}}
\newcommand{\brem}{\begin{rem}}
\newcommand{\erem}{\end{rem}}
\newcommand{\bex}{\begin{ex}}
\newcommand{\eex}{\end{ex}}
\newcommand{\bcon}{\begin{con}}
\newcommand{\econ}{\end{con}}

\usepackage[textheight=220mm,textwidth=150mm]{geometry}

\title{Enriques involutions and Brauer classes}
\author{Alexei N. Skorobogatov and Domenico Valloni}
\date{\today}
\begin{document}
\maketitle

\begin{abstract}
We prove that every element of order 2 in the Brauer group of 
a complex Kummer surface $X$ descends to an Enriques quotient of $X$.
In `generic' cases this gives a bijection between the set ${\mathcal Enr}(X)$ of 
Enriques quotients of $X$ up to isomorphism and the set of Brauer classes of $X$ of order 2.
For some K3 surfaces of Picard rank $20$ we prove that the fibres of ${\mathcal Enr}(X)\to\Br(X)[2]$
above the non-zero points have the same cardinality.
\end{abstract}

\section{Introduction}

Let $S$ be a complex Enriques surface and let $\pi\colon X\to S$ be its 
K3 \'etale double cover. 
J.-L. Colliot-Th\'el\`ene asked whether the induced map of Brauer groups 
$\pi^*\colon \Br(S)\simeq\Z/2\to\Br(X)$ is injective or zero. A. Beauville has given a necessary and sufficient
condition for the injectivity of $\pi^*$ \cite[Cor.~5.7]{Bea09} and showed that 
the Enriques surfaces $S$ for which this map is zero form a countable union of hypersurfaces in
the moduli space of Enriques surfaces \cite[Cor.~6.5]{Bea09}. Enriques surfaces with injective $\pi^*$
 are used in explicit constructions of Enriques surfaces over $\Q$
for which the Brauer--Manin obstruction fails to control weak approximation \cite{HS}
and the Hasse principle \cite{VV, 5}. Enriques surfaces over $\Q$ such that the map $\pi^*$ is zero
have been constructed in \cite{HSch11,GSch12}.

From a different perspective, one can start with a K3 surface $X$ and consider
the set $\mathcal F(X)\subset\Aut(X)$ of fixed point free involutions $\si\colon X\to X$, which are precisely
the involutions such that
the quotient $X/\si$ is an Enriques surface. In this paper we are interested in the map
$$\Phi_X\colon{\mathcal F}(X) \lra \Br(X)[2],$$
which sends $\si\in\mathcal F(X)$ to $\pi^*(b_S)$, where
$\pi\colon X\to X/\si=S$ is the quotient morphism, and
$b_S$ is the unique non-zero element of $\Br(S)$. Combination of results of Beauville
and of Keum and Ohashi shows that $\Im(\Phi_X)$ depends only on the transcendental lattice $T(X)$ of $X$,
see Corollary \ref{U}.

Let ${\mathcal Enr}(X)$ be the set of Enriques quotients of $X$, considered up to isomorphism
of varieties. Equivalently, ${\mathcal Enr}(X)$ is the set of conjugacy classes of $\Aut(X)$
contained in $\mathcal F(X)$, see \cite[Prop.~2.1]{Oha07}.
H. Ohashi proved that the set ${\mathcal Enr}(X)$ is always finite \cite[Cor.~0.4]{Oha07}
although its size is not bounded \cite[Thm.~0.1]{Oha07}. The map $\Phi_X$ is $\Aut(X)$-equivariant,
where $\Aut(X)$ acts on $\mathcal F(X)$ by conjugation, so $\Phi_X$ descends to a map
$$\varphi_X\colon{\mathcal Enr}(X) \lra \Br(X)[2]/\Aut(X).$$
The action of $\Aut(X)$ on $\Br(X)[2]$ factors through the action of the group 
of Hodge isometries of the integral Hodge structure on $T(X)$, so when $\Aut_{\rm Hdg}(T(X))=\{\pm 1\}$ 
the action of $\Aut(X)$ on $\Br(X)[2]$ is trivial. 
In such a `generic' situation, $\varphi_X$ is a map ${\mathcal Enr}(X) \to \Br(X)[2]$.
Note that in this case the set ${\mathcal Enr}(X)$ depends only on the isometry class of the lattice $T(X)$, see the discussion after Theorem \ref{KO}.

Examples show that the set 
${\mathcal Enr}(X)$ can be empty or very large, so in general $\varphi_X$ is neither surjective nor injective.
A very general Enriques surface $S$ (corresponding to the points of the moduli space
outside a countable union of hypersurfaces) is the
unique Enriques quotient of its K3 cover $X$; by Beauville,
in this case $\varphi_X({\mathcal Enr}(X))$ is a certain non-zero element of $\Br(X)[2]$.

The aim of this paper is to clarify the structure of $\Phi_X$ and $\varphi_X$ in some favourable situations.
J.H. Keum proved that every Kummer surface is a double cover of some Enriques surface \cite[Thm.~2]{K}.
His method can be used to prove the following 

\medskip

\noindent{\bf Theorem A}. {\em
Let $X$ be a Kummer surface. Then for every $\alpha\in \Br(X)$ of order $2$ there is an Enriques
quotient $\pi_S\colon X\to S$ such that $\alpha=\pi_S^*(b_S)$.}

\medskip

In other words, for Kummer surfaces the set $\Br(X)[2]\setminus\{0\}$ is contained
in the image of $\Phi_X$. As a kind of partial converse,
in Corollary \ref{coco} we show that if $X$ is a K3 surface such that the abelian group $\Br(X)[2]$ 
is generated by the image of $\Phi_X$, then the transcendental lattice of $X$ is divisible by 2
as an even lattice. We do not know if there exist Kummer surfaces
such that $\Phi_X^{-1}(0)$ is non-empty. At the end of Section \ref{pre} we give examples
of non-Kummer K3 surfaces such that $\Im(\Phi_X)=\{0\}$.

In two `generic' cases Ohashi classified all Enriques quotients of a given K3 surface. 
Combining Theorem A with his results \cite[Thm.~4.1]{Oha07}, \cite[Thm.~1.1]{Oha09} we obtain

\medskip

\noindent{\bf Corollary B}. {\em
Let $X$ be the Kummer surface attached to any of the following abelian surfaces:

\smallskip

{\rm (i)} a product of two non-isogenous elliptic curves;

{\rm (ii)} the Jacobian $J$ of a curve of genus $2$ such that $\NS(J)\cong\Z$.

\smallskip

\noindent Then $\varphi_X$ is a bijection ${\mathcal Enr}(X)\tilde\lra \Br(X)[2]\setminus\{0\}$.}
\medskip

For some K3 surfaces of maximal Picard rank the following result gives information about the fibres
of $\varphi_X$. Its proof uses a certain Galois action on $\Br(X)[2]$ constructed
by the second named author in \cite{V}. 

\medskip

\noindent{\bf Theorem C}. {\em
Let $X$ be a K3 surface of Picard rank $20$.
Let $E=\Q(\sqrt{-d})$, where $d$ is the discriminant of the transcendental lattice $T(X)$.
Assume that $\End_{\mathrm{Hdg}} (T(X))$ is the ring of integers $\O_E\subset E$ and, moreover,
$2$ is inert in $E$ and $E \neq \Q(\sqrt{-3})$.
Then $\Aut_{\mathrm{Hdg}} (T(X))=\{\pm1\}$ and the fibres of 
$\varphi_X\colon{\mathcal Enr}(X) \to \Br(X)[2]$ above the non-zero points
have the same cardinality.}

\medskip

The conditions in Theorem C are easy to check. Let
$$\left(\begin{array}{cc} 2a&b\\b&2c\end{array}\right)$$
be the Gram matrix of $T(X)$, where $a,b,c\in\Z$, so that $-d=b^2-4ac<0$.
Write $-d=f^2D$, where $f\in\Z$ and $D$ is the discriminant of $E$. By \cite[Thm.~3.2]{V} we have
$\End_{\mathrm{Hdg}} (T(X))=\O_E$ if and only if $f={\rm gcd}(a,b,c)$.
Next, $2$ is inert in $E$ if and only if $D\equiv 5 \bmod 8$. If $f$ is odd, 
so that $-d\equiv 5 \bmod 8$, we have 
${\mathcal Enr}(X)=\emptyset$ by \cite{S}, so in this case the fibres of $\varphi_X$ are empty.
Using Theorem A it is easy to see that for each $D\equiv 5 \bmod 8$, $D\neq -3$, there are infinitely many 
pairwise non-isomorphic
K3 surfaces of Picard rank $20$ with complex multiplication by $\O_{\Q(\sqrt{D})}$ such that the fibres of $\varphi_X$ 
 above the non-zero points of $\Br(X)[2]$ have the same {\em positive} number of elements.

\medskip

It would be interesting to describe the K3 surfaces $X$ such that $\Phi_X$ is 
surjective onto $\Br(X)[2]$ or onto $\Br(X)[2]\setminus\{0\}$. 
In this direction we have the following result, whose proof uses
Nikulin's theory of lattices \cite{N} and surjectivity of the period map for K3 surfaces. 

\medskip  

\noindent{\bf Theorem D}. {\em
Let $X$ be a K3 surface such that $\rk(\NS(X))\geq 12$. Then there exist infinitely many K3 surfaces $Y$ such that
\begin{enumerate}
    \item $T(X)_\Q \cong T(Y)_\Q$ as polarised Hodge structures;
    \item the discriminants of $T(Y)$ are pairwise different;
    \item there is a natural isomorphism $\Br(X)[2]\cong\Br(Y)[2]$
under which $\Im(\Phi_X)\setminus\{0\}$ is a subset of $\Im(\Phi_Y)\setminus\{0\}$.
\end{enumerate}
}

\medskip

We recall results of Beauville, Keum and Ohashi, and then prove some useful lemmas in Section \ref{pre}.
Theorem A and Corollary B are proved in Section \ref{kumm}, Theorem C is proved in
Section \ref{sing}, and Theorem D in Section \ref{con}.

This note was written while the first named author was a fellow of FRIAS, University of Freiburg.
He is greatful to FRIAS for hospitality and excellent working conditions.
The second named author was funded by the European Research Council under EU Horizon 2020 
research and innovation programme grant agreement 948066.
 
\section{Lattices and the topology of Enriques quotients} \label{pre}

A {\em lattice} $L$ is a free finitely generated abelian group with a non-degenerate integral
symmetric bilinear form. Write $L(2)$ for the same group with the form $2(x.y)$.

For a lattice $L$ we denote by $A_L=L^*/L$ the discriminant group of $L$.
If $L$ is even, then $q_L\colon A_L\to\Q/2\Z$ is the associated quadratic form.

If $L\subset M$ are lattices, we denote by $L^\perp_M$ the orthogonal complement to $L$ in $M$.
It is clear that $L^\perp_M$ is a primitive sublattice of $M$.

\medskip

Let $\UU$ be the hyperbolic plane. Write $\UU=\Z e\oplus\Z f$, where $(e^2)=(f^2)=0$, $(e.f)=1$.
We denote by $\EE_8$ the negative-definite, even, unimodular lattice of the root system $\EE_8$.
Write 
$$\Lambda=\EE_8^{\oplus 2}\oplus \UU^{\oplus 3}, \quad\quad
M=\UU(2)\oplus \EE_8(2), \quad\quad N=\UU\oplus \UU(2)\oplus \EE_8(2).$$
Here $\Lambda$ is the K3 lattice. Let $\iota\colon\Lambda\to\Lambda$
be the involution permuting two copies of 
$\EE_8\oplus \UU$, and acting as $-1$ on the third copy of $\UU$.
Then $\Lambda^+\cong M$ and $\Lambda^-\cong N$, where $\Lambda^{\pm}$ 
is the $\pm 1$-eigenspace of $\iota$.
For any Enriques quotient $\pi_S\colon X\to S=X/\sigma$, the induced map
$$\pi_S\colon \H^2(S,\Z)/_{\rm tors}\,\lra \H^2(X,\Z)$$ 
can be identified with the composition
$$\H^2(S,\Z)/_{\rm tors}\simeq\UU\oplus \EE_8\stackrel{\rm diag}\lra (\UU\oplus \EE_8)^{\oplus 2}\subset \Lambda\simeq\H^2(X,\Z).$$
Then the fixed point free involution $\sigma\colon X\to X$ induces the involution $\iota$ on $\Lambda$.

The lattice $N$ has a canonical character $N\to\Z/2$ which will play a crucial role in what follows.

\ble \label{e}
The homomorphism $\e\colon N\to\Z/2$ given by
$\e(x)=\big(x.(e+f)\big)$, where $e$ and $f$ are standard generators of $\UU\subset N$, does not depend on
the embedding of lattices $\UU\hookrightarrow N$. Hence $\alpha^*(\e)=\e$ for any $\alpha\in\Aut(N)$.
\ele
{\em Proof.} Let $e',f'$ be standard generators of $\UU$ embedded in $N$.
Write $e'=ae+bf+u$, $f'=ce+df+w$, where $a,b,c,d\in\Z$ and $u,w\in \UU(2)\oplus \EE_8(2)$.
We have $2ab+(u^2)=2cd+(w^2)=0$ and $ad+bc+(u.w)=1$. Since $(u^2)$ and $(w^2)$
are divisible by 4, and $(u.w)$ is even, we see
that $ab$ is even, $cd$ is even, and $ad+bc$ is odd. It follows that either $a,d$ are odd and $b,c$ are even,
or $a,d$ are even and $b,c$ are odd. In both cases $e'+f'$ equals $e+f$ modulo 
$2U\oplus \UU(2)\oplus \EE_8(2)$, hence the result. $\Box$

\ble \label{ee}
If $x\in N$ is such that $(x^2)\equiv 2\bmod 4$, then $\e(x)=0$.
\ele
{\em Proof.} Write $x=ae+bf+u$ where $a,b\in\Z$ and $u\in \UU(2)\oplus \EE_8(2)$.
Then $a$ and $b$ are both odd, hence $\e(x)\equiv a+b\equiv 0\bmod 2$. $\Box$

\ble \label{f}
Let $L$ be a sublattice of $N$. If the restriction of $\e\colon N\to\Z/2$ to $L$ is non-zero, then
$L^\perp_N=L'(2)$ for some even lattice $L'$.
\ele
{\em Proof.} Suppose $\e(x)\neq 0$ for some $x\in L$.
Writing $x=ae+bf+u$ where $a,b\in\Z$ and $u\in \UU(2)\oplus \EE_8(2)$, we see
that $a$ and $b$ have opposite pairity. If $y=ce+df+w \in L^\perp_N$, 
 where $c,d\in\Z$ and $w\in \UU(2)\oplus \EE_8(2)$, then
$ad+bc$ is even, which implies that either $c$ or $d$ is even. 
Then $(y^2)=2cd +(w^2)$ is divisible by $4$, hence $L^\perp_N=L'(2)$ for some even lattice $L'$.
$\Box$

\medskip

The importance of the character $\e\colon N\to\Z/2$ has been revealed by Beauville.
Namely, let $\pi_S\colon X\to S$ be an Enriques quotient of a K3 surface $X$.
Let $T=T_X\subset \Lambda$ be the transcendental lattice of $X$. Recall the canonical isomorphism
$$\Br(X)\cong\Hom(T(X),\Q/\Z).$$
The involution $\sigma^*=\iota$ acts on $T(X)$ as $-1$, so $T(X)\subset N$.

\bthe[Beauville] \label{B}
Let $\pi_S\colon X\to S$ be an Enriques quotient of a K3 surface $X$. Then
$\pi_S(b_S)\in\Br(X)[2]$ is the restriction of $\e\colon N\to\Z/2$ to $T(X)$.
\ethe
{\em Proof.} See \cite[Prop.~3.4 and 5.3]{Bea09}. $\Box$

\medskip

An embedding $T(X)\subset N$ coming from an Enriques quotient of $X$ is clearly primitive.
The orthogonal complement $T(X)^\perp_N\subset N$
contains no $(-2)$-elements $x$, because
by Riemann--Roch either $x$ or $-x$ is effective, but $\sigma^*$ preserves effectivity.
In fact, these are the only conditions.
Horikawa's theorem on the surjectivity of the period map for Enriques surfaces \cite{Hor78} leads to the following
result. See \cite[Thm.~1]{K}, which was extended in \cite[Prop.~2.1]{Oha09}.

\bthe[Keum, Ohashi] \label{KO}
Let $X$ be a K3 surface. Associating to an Enriques quotient of $X$ a primitive embedding
$T(X)\subset N$ defines
a bijection between ${\mathcal Enr}(X)$ and the set of primitive embeddings
of $T(X)$ into $N$ without $(-2)$-elements in
the orthogonal complement, considered up to a Hodge isometry of $N$.
\ethe
Here $N$ is given the Hodge structure such that $T(X)^\perp_N\subset N^{1,1}$.

\medskip

If $\Aut_{\rm Hdg}(T(X))=\{\pm 1\}$ (which holds, for example, when the Picard number of $X$ is odd), the set
$\mathcal Enr(X)$ depends only on the lattice $T(X)$. Then Theorem \ref{KO} says that the elements of 
$\mathcal Enr(X)$ bijectively correspond to the primitive embeddings 
$T(X)\hookrightarrow N$ without $(-2)$-elements in
the orthogonal complement, considered up to an {\em isometry} of $N$.

\bco \label{U}
For any K3 surface $X$ the following statements hold.

\smallskip

{\rm (i)} $\Im(\Phi_X)\setminus\{0\}$ is the set of non-zero $\alpha\in\Br(X)[2]\cong\Hom(T(X),\Z/2)$,
for which there exists a primitive embedding $i\colon T(X)\hookrightarrow N$ 
such that $\alpha=i^*(\e)$.

{\rm (ii)} $0\in\Im(\Phi_X)$ if and only if
there exists a primitive embedding $i\colon T(X)\hookrightarrow N$ without $(-2)$-elements in
the orthogonal complement such that $i^*(\e)=0$.

{\rm (iii)} If $x\in T(X)$ is such that $(x^2)\equiv 2\bmod 4$, then $\alpha(x)=0$
for any $\alpha\in\Im(\Phi_X)$.
\eco
{\em Proof.} Parts (i) and (ii) formally follow from Theorems \ref{B}, \ref{KO} and Lemma \ref{f}. 
Part (iii) follows from Lemma \ref{ee}. $\Box$

\bco \label{coco}
If $X$ is a K3 surface such that the abelian group $\Br(X)[2]$ is generated by
the image of $\Phi_X$, then there is an even lattice $T'$ such that $T(X)\cong T'(2)$.
\eco 
{\em Proof.} It is enough to show that for every $x\in T(X)$ we have $(x^2)\equiv 0\bmod 4$.
Suppose that there is an element $y\in T(X)$ such that $(y^2)\equiv 2\bmod 4$.
Then $y$ is not divisible by 2 in $T(X)$. By Corollary \ref{U} (iii) the non-zero class of $y$
in $T(X)/2T(X)$ is in the kernel of every $\alpha\in\Im(\Phi_X)$. 
Thus $\Im(\Phi_X)$ is contained in a proper subgroup of $\Br(X)[2]$. $\Box$

\bco \label{11jan}
Let $X$ be a K3 surface such that $T(X)$ has a basis $e_1,\ldots,e_n$ with
$(e_i^2)\equiv 2 \bmod 4$ for $i=1,\ldots,n$. Then either $\mathcal Enr(X)=\emptyset$ or
$\Im(\Phi_X)=\{0\}$.
\eco
{\em Proof.} Suppose that a non-zero $\alpha\in \Hom(T(X),\Z/2)$ is in the image of
$\Phi_X$. By Theorem \ref{B} there is a primitive
embedding $i\colon T(X)\to N$ such that $i^*(\e)=\alpha$. By Lemma \ref{ee}
we have $\alpha(e_i)=0$ for $i=1,\ldots,n$, hence $\alpha(T(X))=0$ which is a contradiction. $\Box$

\medskip

This can be used to give examples of K3 surfaces $X$ such that $\Im(\Phi_X)=\{0\}$.
For example, one can take the K3 surface $X$ of Picard rank 20 with transcendental lattice
$$\left(\begin{array}{cc}2&0\\0&2c\end{array}\right)$$ with $c=3, 5, 7$.
Indeed, by \cite[Table 3.1]{SV20} in these cases we have $|\mathcal Enr(X)|=1$.

\section{Kummer surfaces} \label{kumm}

\subsubsection*{Proof of Theorem A}

By Corollary \ref{U} (i) it is enough to construct, for any non-zero $\alpha\in\Hom(T(X),\Z/2)$,
a primitive embedding $i\colon T(X)\hookrightarrow N=\UU\oplus\UU(2)\oplus \EE_8(2)$
such that $\e(x)=\alpha(x)$ for any $x\in T(X)$. We use Morrison's classification of 
transcendental lattices of Kummer surfaces, see \cite[Ch.~14, Cor.~3.2]{H}.
For each of them
Keum \cite[pp.~106-108]{K} constructed a primitive embedding into $N$;
we follow the same strategy to construct all $2^n-1$ embeddings, where $n=\rk(T(X))$.
We keep the notation of \cite{K}, in particular,
$e,f$ is a standard basis of $\UU$ and $h,k$ is a standard basis of $\UU(2)$.
We denote by $\rho$ the Picard rank of $X$.

In the proof below we shall use the following particular case of a result of Nikulin.

\ble
Any even negative definite lattice of rank at most $4$ has a primitive embedding in $\EE_8$.
\ele
{\em Proof.} This follows from \cite[Thm.~1.12.4]{N} using the fact that $\EE_8$ is a unique
even unimodular negative definite lattice of rank 8. $\Box$

\subsubsection*{$\rho=20$}

In this case the lattice $T=\Z x\oplus\Z y$ is positive definite with Gram matrix
$$\left(\begin{array}{cc}4a&2b\\2b&4c
\end{array}\right)$$
where $a,b,c\in\Z$. The three primitive 
embeddings can be given by sending $x,y$ to the following two elements of $N$:
$$(e+2af,2bf+h+ck), \quad (2bf+h+ak,e+2cf), \quad (e+2af, e+(2b-2a)f+h+(c-b+a)k).$$

\subsubsection*{$\rho=19$}

Now $T$ has signature $(2,1)$. We can choose an integral basis $x,y,t$ of $T$ so that
the Gram matrix is
$$\left(\begin{array}{ccc}4a&2d&2l\\2d&4b&2m\\2l&2m&4c
\end{array}\right)$$
where $a,b,c,d,l,m\in\Z$ and $a,b,c<0$.
The embeddings we need to construct are numbered
by the non-zero vectors $(v_1,v_2,v_3)\in(\F_2)^3$ given by evaluating $\e$ 
on the images of $x,y,t$ in this order.
By symmetry it is enough to construct embeddings
labelled $(1,0,0)$, $(1,1,0)$, and $(1,1,1)$. The first two
can be given by sending $x,y,t$ to the following three elements of $N$,
where $w$ is a primitive element of $\EE_8(2)$ such that $(w^2)=4c$:
$$(e+2af,2df+h+bk,2lf+mk+w); $$
$$(e+2af,e+(2d-2a)f+h+(b-d+a)k,2lf+(m-l)k+w).$$
Next we deal with $(1,1,1)$. Without loss of generality we can assume $m>0$. Take
$$(e+k+ah,e+2mf+(d-m)h+w', e+lh+w),$$
where $\Z w'\oplus\Z w$ is a primitive sublattice of $\EE_8(2)$
such that $(w'^2)=4b-4m<0$, $(w^2)=4c<0$, $(w'.w)=0$.

\subsubsection*{$\rho=18$}

Here the lattice $T$ is the orthogonal direct sum of $\Z x\oplus\Z y$ with signature $(1,1)$
and Gram matrix 
$$\left(\begin{array}{cc}4a&2b\\2b&4c\end{array}\right)$$
and $\UU(2)=\Z r\oplus\Z s$. Without loss of generality we assume that $a,c<0$ and $b>0$.
Let $w$ and $u$ be primitive vectors of $\EE_8(2)$ such that $(w^2)=4c<0$
and $(u^2)=4(a-b+c)<0$. We label the embeddings in the same way as above.

Let us first construct primitive embeddings with labels $(1,0,0,0)$ and $(1,1,0,0)$
by taking the direct sum of a primitive embedding $\Z x\oplus\Z y$ into $\UU\oplus\EE_8(2)$ and 
the identity embedding $\UU(2)\tilde\lra\UU(2)$. We send $x,y$ to
$$(e+2af,2bf+w), \quad (e+2af,e+(2b-2a)f+u).$$
The embedding labelled $(0,1,0,0)$ is obtained from $(1,0,0,0)$ by exchanging the roles of $x$ and $y$.

The embedding with label $(1,0,1,0)$ can be obtained by sending $x,y,r,s$ to
$$(e+2af-ak,2bf-bk+w,e+h,k).$$
For $(0,0,1,0)$ we take $(h+w_1,bk+w_2,e,2e+2f+w_3)$, where 
$\Z w_1\oplus\Z w_2\oplus \Z w_3$ is a primitive sublattice of $\EE_8(2)$ with diagonal Gram matrix
such that $(w_1^2)=4a<0$, $(w_2^2)=4c<0$, $(w_3^2)=-8$.

For $(0,0,1,1)$ we take $(h+w_1,bk+w_2,e,e+2f+w_3)$, where 
$\Z w_1\oplus\Z w_2\oplus \Z w_3$ is a primitive sublattice of $\EE_8(2)$ with diagonal Gram matrix
such that $(w_1^2)=4a<0$, $(w_2^2)=4c<0$, $(w_3^2)=-4$.

For $(1,1,1,0)$ we take $(e+2af-ak,e+(2b-2a)f+(a-b)k+u,e+h,k)$.

For $(1,0,1,1)$ we take $(e+2af-ak,2bf-bk+w,e+h,e+k+h+w')$, where
$\Z w\oplus\Z w'$ is a primitive sublattice of $\EE_8(2)$ 
such that $(w^2)=4c<0$, $(w'^2)=-4$, $(w.w')=0$.

For $(1,1,1,1)$ we take $(e+2af-ak,e+(2b-2a)f+(a-b)k+u,e+h,e+k+h+w')$, where
$\Z u\oplus\Z w'$ is a primitive sublattice of $\EE_8(2)$ 
such that $(u^2)=4(a-b+c)<0$, $(w'^2)=-4$, $(u.w')=0$.

The remaining labels can be obtained by swapping $x$ and $y$ as well as $r$ and $s$.

\subsubsection*{$\rho=17$}

Here we have $T=\UU(2)\oplus\UU(2)\oplus(-4m)$, where $m\geq 1$. 
A standard basis is $\{x,y,x',y',t\}$. Up to swapping the two copies of $\UU(2)$ and swapping
the elements of a standard basis of each $\UU(2)$ it is enough to construct the embeddings with
the following labels:
$$(1,0,0,0,0), (1,1,0,0,0), (1,0,0,0,1), (1,1,0,0,1), (0,0,0,0,1), $$
$$(1,1,1,1,0), (1,1,1,1,1), (1,0,1,0,0), (1,1,1,0,0), (1,1,1,0,1), (1,0,1,0,1).$$
The first five embeddings are obtained as direct sums of a primitive embedding of 
$\UU(2)\oplus(-4m)$ into $\UU\oplus\EE_8(2)$ and the identity embedding $\UU(2)\tilde\lra\UU(2)$.
The respective primitive embeddings of $\UU(2)\oplus(-4m)$ into $\UU\oplus\EE_8(2)$
are given by sending $x,y,t$ to the following triples:
$$(e,2e+2f+u_1, v_1), (e,e+2f+u_2,v_2), (e,2e+2f+u_3,e+v_3), (e,e+2f+u_4,e+v_4),$$
Here $\Z u_i\oplus\Z v_i$ is a primitive sublattice of $\EE_8(2)$ such that

$(u_1^2)=-8$, $(v_1^2)=-4m$, $(u_1.v_1)=0$;

$(u_2^2)=-4$, $(v_2^2)=-4m$, $(u_2.v_2)=0$;

$(u_3^2)=-8$, $(v_3^2)=-4m$, $(u_3.v_3)=-2$;

$(u_4^2)=-4$, $(v_4^2)=-4m$, $(u_4.v_4)=-2$.

\noindent The embedding labelled $(0,0,0,0,1)$ can be obtained by sending $x,y,t$ to
$$(2e+2f+w_0,2e+2f+w_1,e+w_2),$$
where $w_0,w_1,w_2$ generate a primitive sublattice of $\EE_8(2)$ with Gram matrix
$$\left(\begin{array}{ccc}-8&-6&-2\\-6&-8&-2\\-2&-2&-4m
\end{array}\right)$$
Indeed, this matrix is negative definite.

To construct the last six embeddings we exhibit the images of $x,y,x',y',t$.
In the case of $(1,1,1,1,0)$ we consider
$$(e,e+2f+k+w_0,e-h,e-h-k+w_1,w_2),$$
where $w_0,w_1,w_2$ generate a primitive sublattice of $\EE_8(2)$ with diagonal Gram matrix
such that $(w_0^2)=(w_1^2)=-4$ and $(w_2^2)=-4m$.

In the case of $(1,1,1,1,1)$ we take
$$(e,e+2f+k+w_0,e-h,e-h-k+w_1,e+w_2),$$
where $w_0,w_1,w_2$ generate a primitive sublattice of $\EE_8(2)$ with the negative definite Gram matrix
$$\left(\begin{array}{ccc}-4&0&-2\\0&-4&0\\-2&0&-4m
\end{array}\right)$$

In the case of $(1,0,1,0,0)$ we take $(e,2f+k,e-h,-k,w)$,
where $w$ is a primitive element of $\EE_8(2)$ with $(w^2)=-4m$.

For $(1,1,1,0,0)$ we take $(e,e+2f+k+u_2,e-h,-k,v_2)$.

For $(1,1,1,0,1)$ we take $(e,e+2f+k+u_4,e-h,-k,e+v_4)$.

For $(1,0,1,0,1)$ we take $(e,2e+2f+k+u_3,e-h,-k,e+v_3)$.

\subsubsection*{Proof of Corollary B}

(i) 
Let $E_1$ and $E_2$ be non-isogenous elliptic curves and let $X=\Kum(E_1\times E_2)$.
By \cite[\S 4]{Oha07} we have
$\Aut_{\rm Hdg}(T(X))=\{\pm 1\}$ and $|\mathcal Enr(X)|=15$.
In this case all the 15 Enriques involutions of $X$ can be described geometrically: these are
Lieberman involutions and Kondo--Mukai involutions. In this case $\rk(T(X))=4$,
so we have $|\Br(X)[2]\setminus\{0\}|=15$.
\smallskip

\noindent (ii) 
Let $C$ be a smooth projective curve of genus 2 such that $\NS({\rm Jac}(C))\cong\Z$.
Let $X=\Kum({\rm Jac}(C))$.
Condition $\Aut_{\rm Hdg}(T(X))=\{\pm 1\}$ is satisfied since the Picard rank of $X$ is odd.
Ohashi shows in \cite{Oha09} that $|\mathcal Enr(X)|=31$
and describes these 31 involutions geometrically.
In this case $\rk(T(X))=5$, so $|\Br(X)[2]\setminus\{0\}|=31$.

\smallskip

Taking into account (i) and (ii), Corollary B follows from Theorem A 
since a surjective map of finite sets of the same cardinality is a bijection. $\Box$

\section{Singular K3 surfaces} \label{sing}

\subsubsection*{K3 surfaces over $\ov \Q$}

For a variety $X$ over $\ov\Q$ and an element $g \in \Gal(\ov\Q/\Q)$ we
define $X^g=X\times_{\ov\Q, g}\ov\Q$. Then we have a morphism $g\colon X\to X^g$ making the
following diagram commutative
$$\xymatrix{ X\ar[r]^g\ar[d]& X^g\ar[d]\\
\Spec(\ov\Q)\ar[r]^{(g^{-1})^*}&\Spec(\ov\Q)}$$
Here the vertical arrows are structure morphisms.
A morphism of $\ov\Q$-varieties $\phi\colon X\to Y$ gives rise to a morphism 
of $\ov\Q$-varieties $\phi^g=g\phi g^{-1}\colon X^g\to Y^g$.

Let $K\subset\ov\Q$ be a subfield and let $G_K=\Gal(\ov\Q/K)$.

\bde \label{field}
Let $X$ be a variety over $\ov\Q$.

{\rm (i)} The field of moduli of $X$ over $K$ is the subfield $K(X)\subset\ov\Q$ fixed by the group
$\{g\in G_K | X\cong X^g \}$. 

{\rm (ii)}
Let $B \subset \Br(X)$ be a finite subgroup.
The field of moduli of the pair $(X, B)$ over $K$ is the subfield $K(X,B)\subset\ov\Q$ fixed
by the group
$$\{g\in G_K| \exists \, \text{an isomorphism} \, f \colon X^g\rightarrow X \,
\text{such that} \, (g^* \circ f^*)|_B= \id_B \}.$$
\ede

Let us fix an embedding $\ov\Q\subset\C$. For a K3 surface $X$ over $\ov\Q$ we write
$T(X)$ for the transcendental lattice of $X_\C$. One has natural isomorphisms
\cite[Prop.~5.2.3 and p.~142]{CTS21}
$$
\Br(X) \cong \Br(X_\C) \cong \Hom(T(X), \Q / \Z).    
$$

\begin{rem} \label{Remark after field of moduli}
{\rm Let $X$ be a K3 surface over $\ov \Q$
of Picard rank at least $12$. According to \cite[Remark 6.1 (2), p.~32]{V}
a Hodge isometry $h \colon T(X^g) \tilde\rightarrow T(X)$ exists if and only if $X \cong X^g$.
It follows that in this case $K(X,B)$ is the fixed field of the group
    $$ \{g\in G_K | \exists \, \text{a Hodge isometry} \, h \colon T(X^g) \rightarrow T(X) \, \text{such that} \,
(g^* \circ h^*)|_B = \id_B \}. $$ }
\end{rem}

For a K3 surface over $\ov\Q$ we have $\Aut(X)=\Aut(X_\C)$, since $Aut_{X/\ov\Q}$ is a discrete group scheme. Hence the set of conjugacy classes
of fixed point free involutions $\mathcal Enr(X)\subset\Aut(X)$ coincides with $\mathcal Enr(X_\C)$.

\bpr \label{proposition galois action}
Let $X$ be a K3 surface over $\ov \Q$ such that $\Aut_{\rm Hdg}(T(X))=\{\pm 1\}$.
The Galois group $G_{K(X)}$ acts naturally on ${\mathcal Enr}(X)$ and on $\Br(X)[2]$ 
so that the map $\varphi_X\colon{\mathcal Enr}(X) \to \Br(X)[2]$ is $G_{K(X)}$-equivariant.
\epr
{\em Proof.} Write $K := K(X)$. 
We use $\si$ and $\tau$ to denote arbitrary elements of $G_K$.
By Definition \ref{field} (i) we can
find an isomorphism $f_{\si,\tau} \colon X^\si \tilde\lra X^\tau$. 

Let us denote the conjugacy class of $\psi\in \Aut(X)$ by $[\psi]$.

A fixed point free involution $\iota\colon X\to X$ gives rise to a
fixed point free involution $\iota^\si=\si\iota\si^{-1}\colon X^\si\rightarrow X^\si$ for any $\si\in G_K$.
Thus $(\iota^\si)^\tau=\iota^{\tau\si}$ for any $\si, \tau\in G_K$.
 We define an action of $G_K$ on ${\mathcal Enr}(X)$ by making $\si$ send $[\iota]$
to $[f_{1,\si}^{-1}\iota^\si f_{1,\si}]$.
This class depends neither on the choice of $\iota$ in its conjugacy class,
nor on the choice of $f_{1,\si}$. For any $\si, \tau\in G_K$ we have
$$[f_{1,\tau}^{-1}(f_{1,\si}^{-1}\iota^\si f_{1,\si})^\tau f_{1,\tau}]=
[(f_{1,\si}^\tau f_{1,\tau})^{-1}\iota^{\tau\si}(f_{1,\si}^\tau f_{1,\tau})]
=[f_{1,\tau\si}^{-1}\iota^{\tau\si} f_{1,\tau\si}], $$
because $f_{1,\tau\si}$ and $f_{1,\si}^\tau f_{1,\tau}$ are both isomorphisms $X\tilde\lra X^{\tau\si}$,
so replacing one of them by the other does not change the conjugacy class.

Let us now define an action of $G_K$ on $\Br(X)[2]$ by making
$\sigma  \in G_{K}$ act as $f_{1,\si}^*(\si^{-1})^*$ which is 
induced by $\si^{-1}f_{1,\si}\colon X\to X^\si\to X$. This action on $\Br(X)[2]$
does not depend on the choice of $f_{1,\si}$.
Indeed, $f_{1,\si}$ is well defined up to an automorphism of $X$, but the action
of $\Aut(X)$ on $\Br(X)[2]$ factors through the action of $\Aut_{\rm Hdg}(T(X))$.
The latter group is $\{\pm 1\}$ by assumption, so $\Aut(X)$ acts on $\Br(X)[2]$ trivially.
The map $(f_{1,\si})^\tau=\tau f_{1,\si}\tau^{-1}$ is an isomorphism $X^\tau\tilde\lra X^{\tau\si}$,
hence $(f_{1,\si})^\tau f_{1,\tau}$ is an isomorphism $X\to X^{\tau\si}$, so for
the purpose of calculating the induced action of $\Br(X)[2]$ we can replace it with $f_{1,\tau\si}$.
This shows that sending $\sigma  \in G_{K}$ to the map induced on $\Br(X)[2]$ by 
$\si^{-1}f_{1,\si}$ is indeed an action.

We have a commutative diagram
$$\xymatrix{
X\ar[r]^{f_{1,\si}}\ar[d]&X^\si\ar[r]^{\si^{-1}}\ar[d]&X\ar[d]\\
X/(f_{1,\si}^{-1}\iota^\si f_{1,\si})\ar[r]&X^\si/\iota^\si\ar[r]&X/\iota
}$$
where the vertical maps are quotients by the respective fixed point free involutions.
Thus the image of the non-zero element of $\Br(X/\iota)$ in $\Br(X)[2]$ followed by
the action of $\si$ on $\Br(X)[2]$ is the same as the image of the non-zero element of
$\Br(X/(f_{1,\si}^{-1}\iota^\si f_{1,\si}))$ in $\Br(X)[2]$. This proves that 
$\varphi_X$ is $G_{K}$-equivariant. $\Box$

\subsubsection*{Moduli fields of singular K3 surfaces}

Let $X$ is a singular K3 surface, that is, a K3 surface of maximal Picard rank $20$. 
It is well known that every singular K3 surface is defined over $\ov\Q$
and has complex multiplication by the imaginary quadratic field $E = \End_{\mathrm{Hdg}} (T(X)_\Q)$. Assume that $\End_{\mathrm{Hdg}}(T(X))$ is the ring of integers
$\mathcal{O}_E\subset E$. In this situation, the results of \cite{V} give explicit descriptions 
of the moduli fields $E(X)$ and $E(X, \Br(X)[n])$ which we now recall. 

The group $\Br(X)\cong\Hom(T(X),\Q/\Z)$ is naturally an $\mathcal{O}_E$-module.
Let $K_n/E$ be the ray class field of $E$ with modulus $n\mathcal{O}_E$ and let
$\mathrm{Cl}_n(E) \cong \mathrm{Gal}(K_n/E)$. 
The complex conjugation $c$ acts on $\mathrm{Cl}_n(E)$.
Let $\mathrm{Cl}_n(E)^c$ be the $c$-invariant subgroup of $\mathrm{Cl}_n(E)$.
Define
$\widetilde{K}_n \subset K_n$ as the fixed field of $\mathrm{Cl}_n(E)^c$, so that 
$\mathrm{Gal}(\widetilde{K}_n/E) \cong \mathrm{Cl}_n(E)/ \mathrm{Cl}_n(E)^c$. 
Note that $K_1$ is the Hilbert class field of $E$ and $\mathrm{Cl}_1(E) = \mathrm{Cl}(E)$ is the usual class group. The complex conjugation $c$ acts on  $\mathrm{Cl}(E)$ as $-1$. 

\begin{theo} \label{4.6}
Let $X$ be a singular K3 surface. Then $\widetilde{K}_n=E(X, \Br(X)[n])$.
\end{theo}
{\em Proof.} See \cite[Thm.~11.2, Remark 9.2 on p.~41]{V}. $\Box$

\medskip

In particular, we have $\widetilde{K}_1=E(X)$.
If $n$ divides $m$, then $\widetilde{K}_n \subset  \widetilde{K}_m$. 

\subsubsection*{Proof of Theorem C}

The assumptions of Theorem C imply that $\Aut_{\rm Hdg}(T(X))=\O_E^\times=\{\pm 1\}$,
so we can apply Proposition \ref{proposition galois action}. 
Let $\rho$ be the representation of $G_{\widetilde{K}_1}$ in $\Br(X)[2] \cong (\Z/2)^2$
constructed in the proof of Proposition \ref{proposition galois action}. 
It is enough to show that under our assumptions one has $|\rho(G_{\widetilde{K}_1})| = 3$. Then
$G_{\widetilde{K}_1}$ acts transitively on $\Br(X)[2]\setminus\{0\}$, so in view of the 
$G_{\widetilde{K}_1}$-equivariance established in
Proposition \ref{proposition galois action} this will imply Theorem C.
By Theorem \ref{4.6}, we need to prove that $[\widetilde K_2:\widetilde K_1]=3$. 

The following exact sequence describes the ray class group $\mathrm{Cl}_2(E)$:
$$ 0 \rightarrow \frac{\mathcal{O}_E^\times}{\{ x \in \mathcal{O}_E^\times | x \equiv 1 \, \mathrm{mod} \, 2 \}} \rightarrow (\mathcal{O}_E/2)^\times \rightarrow \mathrm{Cl}_2(E) \rightarrow \mathrm{Cl}(E) \rightarrow 0. $$
Under our assumptions we have $\mathcal{O}_E^\times =\{ x \in \mathcal{O}_E^\times  | x \equiv 1 \, \mathrm{mod} \, 2 \} = \{ \pm 1 \}$. 
Since $2$ is inert in $E$, we have $\mathcal{O}_E / 2 \cong \mathbb{F}_4$, thus the sequence above becomes 
$$ 0  \rightarrow \F_4^\times \rightarrow \mathrm{Cl}_2(E) \rightarrow \mathrm{Cl}(E) \rightarrow 0. $$
This is a sequence of $G$-modules, where $G = \{1, c \}$. 
We have $(\F_4^\times)^c =\{1\}$ and $\H^1(G, \mathbb{F}_4^\times) = 0$, hence $\mathrm{Cl}_2(E)^c = \mathrm{Cl}(E)^c$. From this we obtain
the exact sequence $$0 \rightarrow \F_4^\times \rightarrow \mathrm{Gal}(\widetilde K_2/E) \rightarrow \mathrm{Gal}(\widetilde K_1/E) \rightarrow 0.$$
Thus $[\widetilde K_2:\widetilde K_1]=3$, as required. $\Box$

\begin{rem}
{\rm When $2$ is split, a similar argument shows that the $G_{\widetilde{K}_1}$-action on $\Br(X)[2]$ is trivial. }
\end{rem}

\section{Constructing Enriques involutions} \label{con}

For a finite abelian group $G$ we write $\ell(G)$ for the minimal number of generators of $G$.
For a prime $p$ we denote by $G_p$ the $p$-primary subgroup of $G$.
Recall that for a lattice $L$ we write $A_L=L^*/L$ for the discriminant group of $L$. When $L$ is even,
we denote by $q_L\colon A_L\to\Q/2\Z$ the finite quadratic form of $L$.

We need to recall fundamental results of Nikulin about the existence of lattices and their primitive embeddings. 

Let $q \colon A \rightarrow \Q / 2 \Z$ be a finite quadratic form. The signature $\mathrm{sign}(q) \in \Z / 8 \Z$ of $q$ is defined as $(t_+ - t_-) \bmod 8$, where $(t_+, t_-)$ is the signature of any even lattice whose discriminant form is isomorphic to $(A,q)$ (such a lattice always exists and, moreover, this notion is well-defined). One also has 

\begin{equation} \label{sign sum}
	\mathrm{sign}(q \oplus q') = \mathrm{sign}(q) + \mathrm{sign}(q').
\end{equation}

Write $A = \bigoplus_p A_p$, where $p$ ranges over the prime numbers. Then one has quadratic forms $q_p \colon A_p \rightarrow \Q_p / \Z_p$ when $p$ is odd and $q_2 \colon A_2 \rightarrow \Q_2/ 2 \Z_2$ when $p=2$. It is clear that $q$ is the orthogonal direct sum of the forms $q_p$. 

For an odd prime $p$, a finite abelian $p$-group $A_p$, and a quadratic form $q_p\colon A_p\to\Q_p/\Z_p$,
Nikulin showed in \cite[Thm.~1.9.1]{N} that there is a unique $\Z_p$-lattice $K(q_p)$ of rank 
$\ell(A_p)$ whose quadratic form is isomorphic to $q_p$. 

Now let $p=2$. 
Let $q^{(2)}_\theta (2)$ be the discriminant quadratic form of the 1-dimensional $\Z_2$-lattice $(2 \theta)$, 
where $\theta \in \Z_2^\times$.
For a finite abelian 2-group $A_2$ and a quadratic form $q_2\colon A_2\to\Q_2/2\Z_2$
we have the following alternative.
If $q_2$ splits as an orthogonal direct sum $q_2=q^{(2)}_\theta (2) \oplus q'_2$, 
then there are precisely two even $\Z_2$-lattices of rank $\ell(A_2)$ whose quadratic form is isomorphic to $q_2$. If such a splitting of $q_2$ does not exists,
there is a unique $\Z_2$-lattice $K(q_2)$ of rank $\ell(A_2)$ whose quadratic form is isomorphic to $q_2$.
The following result is \cite[Thm.~1.10.1]{N}.

\begin{theo}[Nikulin] \label{Nikulin existence of lattices}
	An even lattice with signature $(t_+, t_-)$ and quadratic form $q\colon A\to\Q/2\Z$ 
exists if and only if the following conditions are satisfied:

\smallskip

\noindent $(1)$ $t_+ - t_- \equiv \mathrm{sign}(q) \mod 8$;

\noindent $(2)$ $t_+, t_- \geq 0$ and $t_+ + t_- \geq \ell(A)$;

\noindent $(3)$ $(-1)^{t_-} |A_p| \equiv \mathrm{discr} K(q_p) \bmod \Z_p^{\times 2}$ for the odd primes $p$ such that $t_+ + t_- = \ell(A_p)$;

\noindent $(4)$ $|A_2| \equiv \pm \mathrm{discr} K(q_2) \bmod \Z_2^{\times 2}$ if $t_+ + t_- = \ell(A_2)$ and $q_2 \neq q_{\theta}^{(2)}(2) \oplus q_2'$ for any $\theta$ and $q_2'$. 
\end{theo}

The following result is a consequence of \cite[Prop.~1.15.1]{N} where we took into account that $N$ is the
unique lattice of signature $(2,10)$ whose quadratic form is isomorphic to $q_N$, see \cite[Cor.~1.13.4]{N}.

\bthe[Nikulin] \label{mrak}
Let $L$ be an even lattice with signature $(2_+, k_-)$ and quadratic form $q_L\colon A_L\to\Q/2\Z$.
The existence of a primitive embedding $L \hookrightarrow N$ is equivalent to the existence
of the following data:
\begin{itemize}
	\item subgroups $H_L \subset A_L$ and $H_N \subset A_N$;
	\item an isomorphism of finite quadratic forms 
$\gamma \colon (H_L,q_L|_{H_L}) \xrightarrow{\sim} (H_N,q_N|_{H_N})$;
	\item an even negative definite lattice $K$ of rank $10-k$;
	\item an isomorphism $\delta$ of finite quadratic forms $-q_K$
and the restriction of $q_L\oplus -q_N$ to $\Gamma_{\gamma}^\perp / \Gamma_{\gamma}$,
where the isotropic subgroup $\Gamma_{\gamma}\subset A_L\oplus A_N$ 
is the graph of $\gamma$ in $H_L\oplus H_N\subset A_L\oplus A_N$.
\end{itemize}
Moreover, if $i \colon L \hookrightarrow N$ is a primitive embedding associated to $(H_L, H_N, \gamma, K, \delta)$, then $K \cong i(L)^{\perp}$. 
\ethe

\begin{rem} \label{remark after Nikulin characterization of embeddings}
{\rm (1) Theorem \ref{mrak} remains true when $\Z$-lattices are replaced by $\Z_p$-lattices, 
for any prime $p$. 

(2) If $f \colon \widetilde{K} \rightarrow K$ is an isometry and $\bar{f} \colon A_{\widetilde{K}} \rightarrow A_K$ is the induced isomorphism, then the primitive embeddings $L\hookrightarrow N$
associated to
$(H_L, H_N, \gamma, K, \delta)$ and to $(H_L, H_N, \gamma, \widetilde{K},\delta \circ \bar{f})$ 
are isomorphic.}
\end{rem}

\begin{defi}
Let $L$ be a lattice.
We say that a sublattice $L' \subset L$ of finite index satisfies condition $(*)$ if 
$$\mathrm{gcd}( 2\mathrm{discr}(L), [L \colon L']) = 1$$ and for each prime $p$ 
not dividing $ 2\mathrm{discr}(L)$ we have $\ell(A_{L',p}) < 12 - \rk(L')$. 
\end{defi}


\bpr \label{inf}
Any lattice of positive rank contains infinitely many distinct sublattices $L'\subset L$
satisfying condition $(*)$.
\epr
{\em Proof.} Let $p$ be any odd prime not dividing $\mathrm{discr}(L)$. 
As is well known, see e.g. \cite[Cor.~1.9.3]{N}, the unimodular $p$-adic lattice $L\otimes\Z_p$
has an orthogonal $\Z_p$-basis $v_1,\ldots,v_n$ such that $(v_i^2)\in\Z_p^\times$ for $i=1,\ldots, n$.
The images of $v_1,\ldots,v_n$ in $(L\otimes\Z_p)/p\cong L/p$
form a basis of this $\F_p$-vector space.
Let $L'\subset L$ be the inverse image of the hyperplane spanned by the images of $v_2,\ldots,v_n$.
Thus $[L:L']=p$, so that $\discr(L')=p^2\discr(L)$. Since $p$ does not divide $\discr(L)$,
we have a canonical isomorphism $A_{L'}\cong A_L\oplus A_{L',p}$. It is enough to check that
$\ell(A_{L',p})=1$, which says that $A_{L',p}$ is cyclic. It is clear that $A_{L',p}\cong A_{L'\otimes\Z_p}$,
so it is enough to prove that $\Hom_{\Z_p}(L'\otimes\Z_p,\Z_p)/(L'\otimes\Z_p)\cong\Z/p^2$.
The $\Z_p$-module $L'\otimes\Z_p$ is freely generated by $pv_1,v_2,\ldots,v_n$, hence
the $\Z_p$-module $\Hom_{\Z_p}(L'\otimes\Z_p,\Z_p)\subset L'\otimes\Q_p$ is freely generated by
$p^{-1}v_1,v_2,\ldots,v_n$, which implies the result. $\Box$

\medskip

Condition $(*)$ says that $[L \colon L']$ is odd, hence the inclusion $L' \subset L$ induces a natural 
isomorphism
\begin{equation} \label{identification Brauers}
	\Hom(L', \Z / 2 \Z) \cong \Hom(L, \Z / 2 \Z)
\end{equation}

Recall that for a primitive embedding $i\colon L\hookrightarrow N$ we denote by $i^*(\e)$ the precomposition of
the character $\e\colon N\to\Z/2$ with $i$.

\begin{theo} \label{main theorem overlattices}
	Let $L' \subset L$ be an inclusion of even lattices of signature $(2_+, k_-)$,
where $k\geq 0$. Then we have the following statements.

\smallskip
	
\noindent{\rm (a)} If $L' \subset L$ satisfies condition $(*)$, then for any primitive embedding
$i\colon L\hookrightarrow N$ with $i^*(\e)\neq 0$ there exists
a primitive embedding $i'\colon L'\hookrightarrow N$ such that $i'^*(\e)=i^*(\e)$
under the identification \eqref{identification Brauers}.

\smallskip

\noindent{\rm (b)} If $[L \colon L']$ is odd, then for any primitive embedding 
$i'\colon L'\hookrightarrow N$ with $i'^*(\e)\neq 0$ there exists a primitive embedding
$i\colon L\hookrightarrow N$ such that $i'^*(\e)=i^*(\e)$
under the identification \eqref{identification Brauers}.
\end{theo}
{\em Proof.} (a) Let $i\colon L\hookrightarrow N$ be a primitive embedding such that $i^*(\e)\neq 0$.
Then $K:=i(L)^\perp_N$ is an even negative definite lattice of rank $10-k$, where $k\leq 10$. By Theorem \ref{mrak}
the embedding $i$ corresponds to some datum $(H_L, H_N, \gamma, K, \delta)$. 

Since $L' \subset L$ satisfies condition $(*)$, the index $[L:L']$ is coprime to $|A_L|$, hence
$A_{L'}$ canonically isomorphic to
$A_L \oplus A_\mathrm{new}$, 
where $|A_\mathrm{new}|=[L:L']^2$. Then $q_{L'}$
is an orthogonal direct sum $q_{L'}\cong q_L\oplus q_{\rm new}$, where $q_{\rm new}$ is a quadratic form
 on $A_\mathrm{new}$.

We claim that there is a negative definite lattice $K'$ of rank $10 - k$
such that $A_{K'}\cong A_K \oplus A_\mathrm{new}$ and $q_{K'}\cong q_K\oplus -q_{\rm new}$.
Since $L'$ is a sublattice of $L$ of finite index and $\rk(K)=10-k$ we have
$$\mathrm{sign}(q_L)\equiv \mathrm{sign}(q_{L'}) \bmod 8, \quad \quad
k-10\equiv\mathrm{sign}(q_K) \bmod 8.$$
Since $q_{L'}\cong q_L \oplus q_\mathrm{new}$, we have that $\mathrm{sign}(q_{L'}) = \mathrm{sign}(q_{L}) + \mathrm{sign}(q_{\mathrm{new}})$ by (\ref{sign sum}). 
Thus $\mathrm{sign}(q_\mathrm{new}) \equiv 0 \bmod 8$, which implies property (1) of Theorem \ref{Nikulin existence of lattices}.  

By condition $(*)$ we know that $|A_\mathrm{new}|$ is odd and coprime to $|A_L|$.
For any odd prime $p$ the $\Z_p$-lattices $L \otimes \Z_p$ and $K \otimes \Z_p$ are orthogonal complements of each other in the unimodular $\Z_p$-lattice $N \otimes \Z_p$, hence $|A_{L,p}|=|A_{K,p}|$.
Thus $|A_K|$ and $|A_\mathrm{new}|$ are coprime. This implies
$$\ell(A_K \oplus A_\mathrm{new}) = \max\{\ell(A_K), \ell(A_\mathrm{new})\}\leq 10-k,$$
since $\ell(A_K)\leq\rk(K)=10-k$ and $\ell(A_\mathrm{new})\leq \ell(A_{L'})< 12-\rk(L)$
by condition $(*)$. Thus property (2) of Theorem \ref{Nikulin existence of lattices} also holds.

We now check properties (3) and (4) taking into account the coprimality of $|A_K|$ and $|A_\mathrm{new}|$.
If $p$ divides $|A_K|$, then (3) or (4) holds because the lattice $K$ exists. 
If $p$ divides $|A_\mathrm{new}|$, then $\ell( A_\mathrm{new}) < \rk(K')$ by 
condition $(*)$, so there is nothing to check.

Theorem \ref{Nikulin existence of lattices} now implies the existence of $K'$ with required properties.

\smallskip

Let us construct a datum defining the desired primitive embedding $L' \hookrightarrow N$. 
Since $2A_N=0$, we have $2H_N=0$ and thus $2H_L=0$, so that $H_L \subset A_{L,2}$. 
In view of the canonical isomorphism $A_{L,2} \cong A_{L',2}$, we can  
keep the same $H_{L'}=H_L$, $H_N$ and $\gamma'=\gamma$ as the first three entries of our datum.  

Recall that $A_{L'}\cong A_L \oplus A_{\mathrm{new}}$. We have 
$$\Gamma_{\gamma'} = \Gamma_\gamma \oplus 0 \subset 
\Gamma_{\gamma'}^\perp = \Gamma_\gamma^\perp \oplus A_\mathrm{new} \subset (A_L \oplus A_N) \oplus A_\mathrm{new},$$
hence $\Gamma_{\gamma'}^\perp / \Gamma_{\gamma'} =
\Gamma_{\gamma}^\perp / \Gamma_{\gamma}\oplus A_\mathrm{new} \cong A_K \oplus A_\mathrm{new}$. 
The restriction of $q_{L'}\oplus -q_N\cong(q_L\oplus -q_N)\oplus q_{\rm new}$ to $\Gamma_{\gamma'}^\perp / \Gamma_{\gamma'}$ is isomorphic
to $-q_K\oplus q_{\rm new}$ via the isomorphism $(\delta, \id)$.

Take a negative definite lattice $K'$ of rank $10 - k$
such that $A_{K'}\cong A_K \oplus A_\mathrm{new}$ and $q_{K'}\cong q_K\oplus -q_{\rm new}$.
Let $\delta'$ be an isomorphism of $-q_{K'}$ with the restriction of $q_{L'}\oplus -q_N$ to 
$\Gamma_{\gamma'}^\perp / \Gamma_{\gamma'}$. 
Let $i'\colon L'\hookrightarrow N$ be a primitive embedding associated to the datum
$(H_{L'}, H_{N}, \gamma', K', \delta')$.

To prove that $i'^*(\e)=i(\e)$ under the natural identification (\ref{identification Brauers}), it is enough to show that the induced embeddings of $\Z_2$-lattices $i_2\colon L\otimes \Z_2\hookrightarrow N\otimes \Z_2$ and 
$i'_2\colon L'\otimes \Z_2\hookrightarrow N\otimes \Z_2$ are isomorphic. 

First we claim that $K \otimes \Z_2$ and $K' \otimes \Z_2$ are isomorphic $\Z_2$-lattices. 
Since $K$ and $K'$ are negative
definite of the same rank, and $|A_{K'}| = |A_K| \cdot |A_\mathrm{new}|$, 
we have $\discr(K')=\discr(K)\cdot |A_{\rm new}|$.
Since $|A_{\rm new}|$ is a square of an odd integer, the even $2$-adic lattices
$K \otimes \Z_2$ and $K' \otimes \Z_2$ have the same rank, discriminant form, and determinant modulo $\Z_2^{\times 2}$. This implies that the $\Z_2$-lattices $K \otimes \Z_2$ and $K' \otimes \Z_2$ are isomorphic 
\cite[Cor.~1.9.3]{N}.

It remains to show that after tensoring with $\Z_2$ the data
$(H_{L}, H_{N}, \gamma, K, \delta)$ and $(H_{L'}, H_{N}, \gamma', K', \delta')$ 
give rise to isomorphic embeddings of 
$L' \otimes \Z_2 \cong L \otimes \Z_2$. The first three entries of each datum are the same. 
By Remark \ref{remark after Nikulin characterization of embeddings} it is enough to find an isometry $f \colon K' \otimes \Z_2 \rightarrow K \otimes \Z_2$ such that $\delta'_2 =\delta_2 \circ \bar{f}$. The existence of 
such an $f$ follows from \cite[Thm.~1.9.5]{N}. This concludes the proof of (a).

\smallskip

(b) Write $A := A_L = A_2 \oplus A_{\mathrm{odd}}$, where $A_2$ is the 2-primary subgroup of $A$. Similarly, write $A' := A_{L'} = A'_2 \oplus A'_{\mathrm{odd}}$. It is clear that $A_2\cong A'_2$.
Then $q_{L'}$ is an orthogonal direct sum of quadratic forms $q_{L,2}$ on $A_2$ and $q_{\rm odd}$
on $A'_{\mathrm{odd}}$.

The overlattice $L$ of $L'$ defines an isotropic subgroup $I \subset A'$, 
where $|I|=[L:L']$, so that $A= I^\perp/I$. Since $[L \colon L']$ is odd by assumption, we have 
$I \subset A'_{\mathrm{odd}}$. Thus $I^\perp = A_2 \oplus I^\perp_{\mathrm{odd}}$, where $I^\perp_\mathrm{odd}=I^\perp\cap A'_\mathrm{odd}$. This shows that $A= A_2 \oplus    I^\perp_{\mathrm{odd}}/I$. 

Let $i'\colon L'\hookrightarrow N$ be a primitive embedding such that $i'^*(\e)\neq 0$.
Then $K':=i(L')^\perp_N$ is an even negative definite lattice of rank $10-k$. Let 
$(H_{L'}, H_N, \gamma', K', \delta')$ be a datum associated to $i'\colon L'\hookrightarrow N$
as in Theorem \ref{mrak}. In particular, $\delta'$ is an isomorphism of $-q_{K'}$
with the restriction of $q_{L'}\oplus -q_N$ to $\Gamma_{\gamma'}^\perp / \Gamma_{\gamma'}$.
Since $2A_N=0$, we have $H_{L'} \subset A'_2=A_2$.
Hence $\Gamma_{\gamma'} \subset A_2\oplus A_N \subset A' \oplus A_N$ and thus $\Gamma_{\gamma'}^\perp = (\Gamma_{\gamma'}^\perp)_2 \oplus A'_{\mathrm{odd}}$, 
where $(\Gamma_{\gamma'}^\perp)_2=\Gamma_{\gamma'}^\perp\cap (A_2 \oplus A_N)$. This shows that $\delta'$ identifies the finite quadratic form $-q_{K'}$ on $A_{K'}$ with the restriction of 
$(q_{L,2}\oplus -q_N)\oplus q_{\rm odd}$ to $\big((\Gamma_{\gamma'}^\perp\big)_2/\Gamma_{\gamma'}) \oplus A'_{\mathrm{odd}}$.

The isotropic subgroup $I \subset A'_\mathrm{odd}$ 
gives rise, via $\delta'$, to an isotropic subgroup in $A_{K'}$, which we still denote by $I$. 
It defines an overlattice $K' \subset K$ with $[K:K']=[L:L']$,
so that $\delta'$ induces an isomorphism $\delta$ of the quadratic form $-q_K$ on $$A_{K} \cong \big( (\Gamma_{\gamma'}^\perp)_2/\Gamma_{\gamma'}\big)  \oplus (I^\perp_{\mathrm{odd}}/I)$$
with the restriction of $(q_{L,2}\oplus -q_N)\oplus q_{\rm odd}$.
Since $K'$ is a sublattice of $K$ of odd index, there is a natural isomorphism
of $\Z_2$-lattices $K' \otimes \Z_2 \cong K \otimes \Z_2$. 

We construct a primitive embedding $i\colon L\hookrightarrow N$ from the datum
$(H_{L}, H_N, \gamma, K, \delta)$, where $H_L=H_{L'}$, $\gamma=\gamma'$, and $K$ and $\delta$ are
defined in the previous paragraph.

To complete the proof of (b) it remains to show that $i$ and $i'$ induce isomorphic embeddings of
$\Z_2$-lattices. This is proved by the same argument as in (a). $\Box$

\begin{cor} \label{cc}
Let $L$ be an even lattice of signature $(2_+, k_-)$ and let $S\subset\Hom(L,\Z/2)$
be a subset such that $0\notin S$. Suppose that for any $\alpha\in S$
there is a primitive embedding $i\colon L\hookrightarrow N$ such that $\alpha=i^*(\e)$.
Then any sublattice of $L$ that satisfies condition $(*)$ and any overlattice of $L$ of odd index
have the same property.
\end{cor}

\subsubsection*{Proof of Theorem D} 

By Proposition \ref{inf} there are infinitely many sublattices $T \subset T(X)$ 
with pairwise different discriminants that satisfy condition $(*)$. 
Endow $T$ with the Hodge structure 
coming from $T(X)$. Since $\rk(T) \leq 10$, by \cite[Thm.~1.14.1]{N} there exists a unique primitive embedding of the lattice $T$ into the K3 lattice $\Lambda$. We equip $\Lambda$ with the Hodge structure induced by the 
Hodge structure on $T$ so that $T^\perp_\Lambda\subset \Lambda^{(1,1)}$.
By the surjectivity of the period map, there is a K3 surface $Y$ together with a Hodge isometry between
$\Lambda$ and $\H^2(Y,\Z)$. The transcendental lattice $T(Y)$ is the orthogonal complement to
$\H^2(Y,\Z)^{(1,1)}$, hence $T(Y)\cong T$. It remains to apply Corollary~\ref{cc}. $\Box$

\bigskip

\noindent Leibniz Universit\"at Hannover,
Riemann Center for Geometry and Physics,
Welfengarten 1,
30167 Hannover, Germany
\smallskip

\noindent {\tt valloni@math.uni-hannover.de}

\bigskip

\noindent Department of Mathematics, South Kensington Campus,
Imperial College London, SW7 2BZ England, U.K. -- and --
Institute for the Information Transmission Problems,
Russian Academy of Sciences, 19 Bolshoi Karetnyi, Moscow 127994
Russia

\smallskip

\noindent {\tt a.skorobogatov@imperial.ac.uk}

\end{document}